\documentclass[11pt,reqno]{article}
\input{epsf.sty}
\usepackage{hyperref}
\usepackage[final]{epsfig}
\usepackage{color}
\usepackage{amsmath, amsthm, amssymb, dsfont, mathrsfs}
\newtheorem {thm}{Theorem}[section]
\newtheorem {prop}[thm]{Proposition}

\def\Cox{\hfill \Box}

\def\N{{\Bbb N}}
\def\TT{{\Bbb T}}
\def\Z{{\Bbb Z}}
\def\R{{\Bbb R}}
\def\P{{\Bbb P}}

\def\E{{\Bbb E}}
\def\one{{\mathds 1}}
\def\Gxi{{\GG(\g[\xi])}}
\def\Gexxi{{{\rm ex}\GG(\g[\xi])}}
\newcommand{\F}{{\mathcal F}}

\newcommand{\C}{{\mathcal C}}

\newcommand{\blue}{\textcolor{blue}}

\def\0{{\bf 0}}

\def\ba{{\backslash}}

\def\oo{\underline{\omega}}

\def\a{\alpha}

\def\b{\beta}

\def\d{\delta}

\def\phi{\varphi}

\def\g{\gamma}

\def\l{\lambda}

\def\s{\sigma}

\def\t{\tau}

\def\x{\xi}

\def\o{\omega}

\def\D{\Delta}

\def\L{\Lambda}

\def\O{\Omega}

\def\T{\T}

\def\AA{{\cal A}}
\def\C{{\cal C}}
\def\DD{{\cal D}}
\def\GG{{\cal G}}
\def\PP{{\cal P}}

\def\FF{{\cal F}}
\def\HH{{\cal H}}
\def\TT{{\cal T}}
\def\EE{{\cal E}}

\def\oo{\underline{\omega}}

\def\OO{\underline{\Omega}}

\def\PP{\underline{P}}

\begin{document}

\title{Gibbs-non Gibbs transitions in different geometries:\\
The Widom-Rowlinson model under stochastic spin-flip dynamics} 

\author{
Christof K\"ulske\footnote{
Ruhr-University of Bochum, Fakult\"at f\"ur Mathematik, Postfach 102148, 44721, Bochum,
Germany
\texttt{Christof.Kuelske@ruhr-uni-bochum.de}}
}

\newcommand{\CC}[1]{{\color{blue} #1}}
\newcommand{\DE}[1]{{\color{red} #1}}
\newcommand{\CK}[1]{{\color{green} #1}}

\maketitle

\begin{abstract}
The Widom-Rowlinson model is an equilibrium model for point particles in 
Euclidean space. It has a repulsive interaction between particles of different colors, 
and shows a phase transition at high intensity. 
Natural versions of the model can moreover be formulated in different geometries:  
in particular as a lattice system or a mean-field system. 
We will discuss recent results on dynamical Gibbs-non Gibbs transitions in this context. 
Main issues will be the possibility or impossibility of an immediate loss of the Gibbs property, 
and of full-measure discontinuities of the time-evolved models. 

(Collaborations with Benedikt Jahnel, Sascha Kissel, Utkir Rozikov)

\end{abstract}

\bigskip 


\section{Introduction}

Recent years have seen a variety of studies of Gibbs-non Gibbs transitions of measures 
which appear as image measures of Gibbs measures, under certain local transformation 
rules.  What is a Gibbs measure? There is a well-defined theory to define Gibbs measures
on lattices, where the probability space is given by the set 
of all functions from lattice sites to a finite alphabet.  The central object is that of a
specification \cite{Ge11,Bo06}.
For other geometries, other but related approaches are adequate, see below. 
The unifying idea  is that Gibbs measures are measures whose conditional probabilites to see 
a single symbol at a given site, are nice (continuous) functions of their conditioning, see below. 
It has been discovered that Gibbs measures under natural 
deterministic or stochastic transformations may lead to non-Gibbsian measures \cite{EnFeSo93}.
More specifically the study of stochastic time evolutions, even very simple ones,
applied to Gibbs measures, has shown very interesting transition phenomena, the most prototypical 
example for this is the Glauber-evolved Ising model in  
\cite{EnFeHoRe02}.   
Indeed, stochastic time evolutions may destroy the Gibbs property of the image measure  
at certain transition times, a phenomenon we call dynamical Gibbs-non Gibbs transitions, see below. 
The purpose of this note is to take the Widom-Rowlinson  model \cite{WiRo70} and variations thereof 
as 
a guiding example, apply an independent symmetric stochastic spin flip 
dynamics to it,
and describe our findings of what may and what may not happen along the time-evolved 
trajectory of measures.  We treat and compare a hard-core version and 
a soft-core version of the model in various geometries, namely in Euclidean 
space, on the lattice, as a mean-field model, and on a regular tree. Our aim here is to provide 
an overview; for detailed statements and proofs we refer to the original articles.

\section{Gibbs on lattice, sequentially Gibbs, marked Gibbs point processes, 
and the Widom-Rowlinson model} 

We start by recalling the notion of an infinite-volume Gibbs measure 
for lattice systems. For the purpose of the discussion of the Widom-Rowlinson model 
and all measures appearing under time-evolution defined below from it, it is sufficient  
to restrict to the local state-space 
$\{-1,0,1\}$ for particles carrying spins plus or minus, and holes. 
Our site space is the lattice $\Z^d$. The space of 
infinite-volume configurations is $\O=\{-1,0,1\}^{\Z^d}$.

\subsection{Specifications and Gibbs measures on the lattice}

The central object in Gibbsian theory on a countable site space 
which defines the model is a {\em specification}. This covers 
both cases of infinite lattices and trees. 
It is a candidate system for conditional probabilities of an 
infinite-volume Gibbs measure $\mu$ (probability measure on $\O$)
to be defined by DLR equations $\mu (\g_\L(f|\cdot))=\mu(f)$. \\
A specification $\g$ is by definition a family of probability kernels
$\g=(\g_\L)_{\L\Subset\Z^d}$, indexed by finite subvolumes $\L$, 
where $\g_\L(d \o|\eta)$ is a probability measure on $\O$, for each fixed 
configuration $\eta$. It must have the following properties. The first is 
the  {\em consistency } which means that 
\begin{equation}\begin{split}\label{einsminus}
\g_\D(\g_\L(d \o|\cdot)|\tilde\o)=\g_\D(d \o|\tilde\o)
\end{split}
 \end{equation}
for all finite volumes
$\L\subset\D\Subset\Z^d$. It is suggested by the tower property of 
conditional expectations.   
 
The second is the  $\F_{\L^c}$-measurability of $\g_\L(f|\cdot)$, for any bounded measurable observable $f$.  
Here the sigma-algebra 
 $\F_{\L^c}$ is generated by the spin-variables outside of the finite volume $\L$. 
 
 The last property is the 
{\em properness} $\g_\L(1_A|\cdot)=1_A$ for $A\in \F_{\L^c}$. It means 
that the randomization of the kernel takes place only inside of $\L$, and 
an event which is determined by what is outside of $\L$ will indeed be determined 
by looking at the boundary condition alone. 

An important additional regularity requirement is {\em quasilocality} of the specification which means 
that  the function $\o\mapsto \g_\L(f|\o)$ should be {\em quasilocal} 
for $f$ quasilocal, and this has to hold for all finite volumes $\L$. 
A quasilocal function is a uniform limit 
of local functions, that is of functions which depend only on finitely many coordinates.  

More specifically a {\em Gibbsian specification} on the infinite-volume state space 
$\O=\{-1,0,1\}^{\Z^d}$ 
for an  {\em interaction potential} $\Phi=(\Phi_A)_{A\Subset \Z^d}$ and a priori 
measure $\alpha \in\mathcal{M}_1(\{-1,0,1\})$ by definition has probability kernels 
\begin{equation}\begin{split}\label{eins}
\g_{\L,\Phi,\alpha}(\o_\L|\o_{\L^{\rm c}}):=\frac{1}{
Z_\L(\o_{\L^{\rm c}})
}e^{-\sum_{A\cap \L\neq \emptyset }\Phi_A(\o)}\prod_{i\in \Lambda}\alpha(\omega_i)
\end{split}
 \end{equation}
where $Z_\L(\o_{\L^{\rm c}})$ is the normalizing partition function. If $\Phi$ is a finite-range potential (meaning that $\Phi_A$ is only nonzero for finitely 
many $A$'s,)  obviously all sums are finite, when we insist that 
$\Phi_A$ takes finite only values. Finiteness of the sums also holds, if $\Phi$ is uniformly absolutely convergent. 
For hard-core models the specification kernels 
acquire an  indicator, see the example \eqref{vier} below.

The first statistical mechanics task in this setup is the following. 
Given a specification $\g=(\g_\L)_{\L\Subset \Z^d}$ in the above sense, find 
the corresponding {\em Gibbs measures} 
\begin{equation}\begin{split}\label{drei}
\GG(\g):=\{\mu\in \mathcal{M}_1(\O) , \mu \g_{\L}=\mu, \text{ for all }\L\Subset \Z^d \}\end{split}
 \end{equation}
In general $\GG(\g)$ may be empty, contain precisely one measure, or more than 
one measures. 
If $|\GG(\g)|>1 $ we say that the specification $\g$ has a {\em phase transition}.  
The Gibbs measures $\GG(\g)$ form a simplex, meaning that each measure has a unique 
decomposition over the extremal elements, called pure states. Pure states can be recovered
as finite-volume limits with fixed boundary conditions. 

Existence and extremal decomposition 
of proper infinite-volume measures 
becomes even more involved for systems with random potentials. In general, 
for systems like spin-glasses, the construction of infinite-volume states by non-random 
sequences of volumes which exhaust the whole lattice is problematic,  
and for such systems the higher-level notion of a {\em metastate} (a measure on infinite-volume 
Gibbs measures) is useful  \cite{NeSt13, Bo06, Ku97, ArDaNeSt10, CoJaKu18}.

\subsection{Hard-core and soft-core Widom-Rowlinson model on lattice and in mean-field}

We will consider here the version of the 
		{\em 	hard-core Widom-Rowlinson model} on $\Z^d$ as in 
		\cite{HiTa04}. 
		It has the a priori measure $\alpha\in\mathcal{M}_1(\{-1,0,1\})$
as its only parameter. Its specification kernels are given by
		\begin{equation}\begin{split}\label{vier}
			\gamma^{hc}_{\Lambda,\alpha}(\omega_\Lambda\vert \omega_{\Lambda^c}) :=\frac{1}{Z^{hc}_\L(\o_{\L^{\rm c}})} I^{hc}_\Lambda(\omega_\Lambda\omega_{\Lambda^c})\textcolor{black}{\prod_{i\in\Lambda}\alpha(\omega_i)},
		\end{split}
 \end{equation}
where the hard-core indicator   $I^{hc}_\Lambda(\omega)= \prod_{i\in\Lambda} I_{(\omega_i\omega_j\neq -1,\,\forall j\sim i)}$ forbids  $+-$ neighbors to occur 
with positive probability. 	
Related hard-core models have been studied on lattices and trees, see for example 
\cite{GaLe71, MaStSu18, Ro13}.  

The 
{\em soft-core Widom-Rowlinson model} on $\Z^d$ 
has an additional repulsion parameter $\beta>0$. 
In the specification kernels, which are by definition given by 
			\begin{align}\label{viera}
			\gamma^{sc}_{\Lambda,\beta,\alpha}(\omega_\Lambda\vert \omega_{\Lambda^c}) :=\frac{1}{Z^{sc}_\L(\o_{\L^{\rm c}})} e^{-\beta \sum_{\{i,j\}\in \mathcal{E}_{\Lambda}^b}I_{(\omega_i\omega_j=-1)} }\textcolor{black}{\prod_{i\in\Lambda}\alpha(\omega_i)}, 
			\end{align}
configurations with $+-$ neighbors are suppressed, but not forbidden. 

These definitions of a specification immediately extend 
to any graph with countably infinite vertex set, 
where each vertex has a finite number of nearest neighbors. 
In particular we may study this model on a regular tree with $k+1$ neighbors, 
see \cite{KiKuRo19}. 

The mean-field formulation is different, as the model is defined 
as a whole sequence of finite-volume Gibbs measures, indexed 
by the system size $N\in\mathbb{N}$. The elements in the sequence for the 
{\em Mean-Field} soft-core Widom-Rowlinson model 
with repulsion parameter $\beta>0$ are the measures 
	\begin{align}
	\mu_{N,\beta,\alpha}(\omega_{[1,N]}):=\frac{1}{Z_{N,\beta,\alpha}} e^{-\frac{\beta}{2N} \sum_{1\leq i,j\leq N} I_{(\omega_i\omega_j =-1)}}\prod_{j=1}^N \alpha(\omega_j)
	\end{align}
for $\omega_{[1,N]}=(\o_i)_{i=1,\dots,N}\in \{-1,0,1\}^N$. For more details see 
\cite{KiKu18}.

\subsection{Sequential Gibbsianness for mean-field (and Kac-models on torus)}

There is an intrinsic formulation of the Gibbs property which is 
suitable also in situations where a finite-volume Hamiltonian  can not 
be read off directly from the explicit definition of the measures.  It focusses 
on conditional probabilities instead, suggested by analogy to the lattice situation
\cite{KuLe07,HoReZu13}. 

Take $(\mu_N)_{N\in \N}$ a sequence of exchangeable probability measures $\mu_N$ 
on the finite-volume state space $\{-1,0,1\}^N$. 
The large $N$-behavior of such a sequence defines our model. 
The model is called {\em sequentially Gibbs} iff the volume-limit 
of the single-spin probabilities in the finite-volume measures  
\begin{equation}\begin{split}\label{fuenfa}\lim_{N\uparrow \infty}\mu_N(d\o_1|\o_{[2,N]})= \g(d \o_1|\nu)\end{split}
 \end{equation} 
 exists whenever the empirical distributions of a configuration
 $(\o_i)_{i=2,3,4,\dots}$ converge, 
\begin{equation}\begin{split}\label{sechs}
\frac{1}{N-1}\sum_{i=2}^N\d_{\o_i}\rightarrow \nu.\end{split}
 \end{equation}
  This has to hold  
 for all {\em limiting 
empirical distributions of conditionings} $\nu\in \mathcal{M}_1(\{-1,0,1\})$. 

If there is some $\nu$ for which it is possible to obtain different limits for 
different boundary conditions $(\o_i)_{i=2,3,\dots}$, and $(\bar\o_i)_{i=2,3,\dots}$
we call this $\nu$ a {\em bad empirical measure}. A model fails to be sequentially 
Gibbs, if there is at least one bad empirical measure. 

As a general consequence of this definition,  the sequential Gibbs property 
implies that the limiting kernel $\nu \mapsto  \g(d \o_1|\nu)$ is continuous. 
In our Widom-Rowlinson case where $\nu$ takes values in a finite-dimensional simplex, 
all topologies are the same, and equal to the Euclidean topology. 
Clearly the mean-field Widom-Rowlinson model defined above, is sequentially 
Gibbs. 

A similar notion of the sequential Gibbs property can be extended to cover 
{\em Kac-models} on the torus, and transformed Kac-models which have the same index set. 
These models are again described by sequences whose asymptotics one wants to capture, 
but they have a spatial structure. 
As in the mean-field models, there is again a single-site limiting kernel, 
however the limiting empirical distribution $\nu$ which appeared as a conditioning 
in the mean-field model 
is replaced by a whole profile of spin densities on the unit torus. For details of these definitions 
and results, see 
\cite{HoReZu13,HeKrKu17}. 

\subsection{Marked Gibbs Point Processes in Euclidean space}

Here the good definition   of  {\em Gibbs measure}   
is in some analogy to the lattice situation
\cite{DeDrGe12,JaKu17,Ru99}. 
We restrict again for the sake of our exposition 
to the specific simple mark space
which covers the Widom-Rowlinson model of point particles 
in Euclidean space, and the time-evolved version we will discuss below. 
In this case the mark space is $\{-1,1\}$. 
It does not contain zero. 
The spatial degrees of freedom are described by the set 
$\O$ of {\em locally finite subsets} of $\R^d$. 
A marked particle configuration is a pair 
$\oo=(\o^{-},\o^{+})$ describing the positions of minus-particles (and plus-particles 
respectively) where each $\o^{-},\o^{+}\in\O$. The configuration
space  of such marked configurations is $\underline{\O}$.  For the measurable 
structure we need the $\s$-algebras 
$\FF,\FF_\L$. These are the  $\s$-algebras for marked particles generated by the counting variables. 
They count the number of plus- and minus particles in Borel sets $A$ in the whole Euclidean space
(or all such sets $A\in \L$ respectively, where $\L$ may be any measurable subset 
of Euclidean space).  

A {\em specification} shall become, as on countable graphs, a 
candidate system for conditional probabilities of Gibbs measure $\mu$ 
to be defined by DLR equations $\mu \g_\L=\mu$ for all measurable bounded 
subsets $\L$ in Euclidean space. 
Hence, one defines a specification to be 
a family of proper probability kernels $\g=(\g_\L)_{\L\Subset\R^d}$ 
with the properties of {\em consistency}, that is  
$\g_\D \g_\L=\g_\D$
for all measurable volumes 
 $\L\subset\D\Subset\R^d$.
 One also needs $\FF_{\L^c}$-measurability of  
 $\g_\L(f|\cdot)$, for any bounded test observable $f$.  
 Properness means here that  $\g_\L(1_A|\cdot)=1_A$ for $A\in \FF_{\L^c}$. 
 
We will further assume 
{\em quasilocality} of the specification. 
This means the compatibility of the kernels $\g_{\L}$ with the {\em local topology} on the space 
of marked point clouds.  In this topology 
convergence for a sequence of marked particle clouds means that the clouds 
must become constant in each bounded volume.

\subsection{Widom-Rowlinson model in Euclidean space}

We assume spatial dimension $d\ge2$, and fix the 
two-color local spin space (mark space) of $\{-,+\}$. The model will 
be obtained as a modification of the base measure $\PP$
by which we denote a 
 two-color homogenous {\em Poisson Point Process} in the infinite volume, with  
intensities $\l_+$ for plus colors and $\l_-$ for minus colors. 

The (hard-core) {\em Widom-Rowlinson specification} is the  Poisson-modification
with the specification kernels 
\begin{equation}\begin{split}\label{sechsa}
\g_\L(d\oo_\L|\oo_{\L^{\rm c}}):=\frac{1}{Z_\L(\oo_{\L^{\rm c}})}\chi(\oo_\L\oo_{\L^{\rm c}})\PP_\L(d\oo_\L)\end{split}
 \end{equation}
where the indicator $\chi$ is one iff the interspecies distance (the distance between points of 
different sign) is bigger or equal than $2a$, 
and $\PP_\L(d\oo_\L)$ denotes the two-color Poisson process in the bounded 
volume $\L$. The picture shows a typical configuration at large $\l_+=\l_-$, in a 
finite volume. 
 \begin{center}
 \includegraphics[scale=0.5]{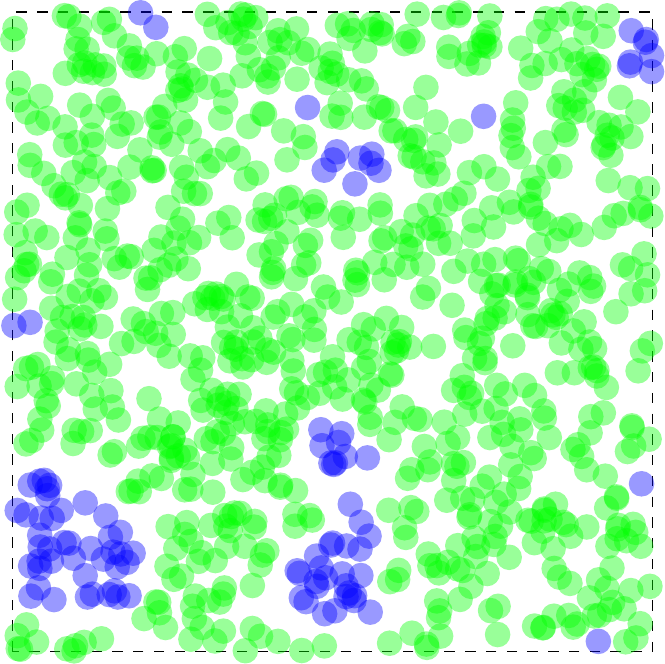} 
 \end{center}
 This picture survives the thermodynamic limit: 
 By results of \cite{Ru71,BrKuLe84, ChChKo95}
 it is known that in 
$d\ge2$, $\l_+=\l_-$ large, the continuum WiRo has a {\em phase transition}, and 
how this is related to percolation of large clusters of overlapping disks. 

In more general models of marked point particles, specifications which are 
Poisson modifications may be obtained in terms 
of exponential factors with finite-volume Hamiltonians which are formed 
with potentials. 
For such potentials one may consider multibody potentials which are known 
from statistical physics, but we may also allow for 
{\em hyperedge potentials } and define 
\begin{equation}\label{PoissonMod}
\begin{split}
\g_\L(d \oo_{\L}|\oo_{\L^c}):=\frac{1}{Z_{\L}(\oo_{\L^c})}e^{
-\sum_{\eta\Subset\o:\, \eta\cap\L\neq\emptyset}\Phi(\eta, \oo)}
\PP_\L(d\oo_{\L}) 
\end{split}
\end{equation}
A hyperedge potential $\Phi(\eta, \oo)$ is by definition 
allowed to depend on the marked point cloud $\oo$ {\em not only } via the particle positions and marks on 
hyperedge $\eta$ (which is just a finite subset 
of points in the cloud), {\em but also  }
on a whole {\em neighborhood } of $\eta$, up to some {\em horizon}. 
This generalization is useful in models of stochastic geometry, 
involving e.g. energies depending on the cells of a 
Voronoi tesselation.  
Hyperedge potentials were successfully used in \cite{DeDrGe12}
where a general existence theory of infinite-volume Gibbs measures is developed. 
We also refer to \cite{JaKu17a} for  
representation theorems. There is it shown how one can go from a measure 
$\mu$ under continuity assumptions of finite-volume conditional probabilities, 
to a hyperedge potential $\Phi$. These theorems
are a generalization of Kozlov-Sullivan theorems 
\cite{Su73,Ko74}
known on the lattice to the continuum, 
and make use in a constructive way 
of the weak nonlocality allowed by the hyperedge potential concept.

\section{Dynamical Gibbs-non Gibbs transitions} 
 
 We now come to time-evolutions. 
 Consider again the Euclidean space Widom-Rowlinson model, and fix 
 some cloud of particles carrying the marks plus or minus. 
We define a continuous-time stochastic dynamics by the following rule.  
Particle locations stay fixed, holes stay fixed. The signs of the particles 
however change stochastically, independently of each other, according to a temporal Poisson process 
for each particle with rate one. In this way, at every particle location, the probability to go 
from $+$ to $-$ in time $t$ is given by 
\begin{equation}\begin{split}\label{sieben}
p_t(+,-)=p_t(-,+)=\frac{1}{2}(1-e^{-2 t}). \end{split}
 \end{equation}
Starting with the same signed particle configuration shown above, 
after a small time $t$, a fraction of the particles has kept their signs, as shown 
in the picture (or flipped back). Of course, there is loss of memory in each fixed bounded volume, 
which is exponentially fast in the time $t$. Interesting things however happen if 
we consider the infinite volume, as we will discuss. 
  \begin{center}
{\hspace{1cm} \includegraphics[scale=0.6]{GnG_PP_WRM_AAP_Version_170410-figure0.pdf} 
\hspace{3cm}
 \includegraphics[scale=0.6]{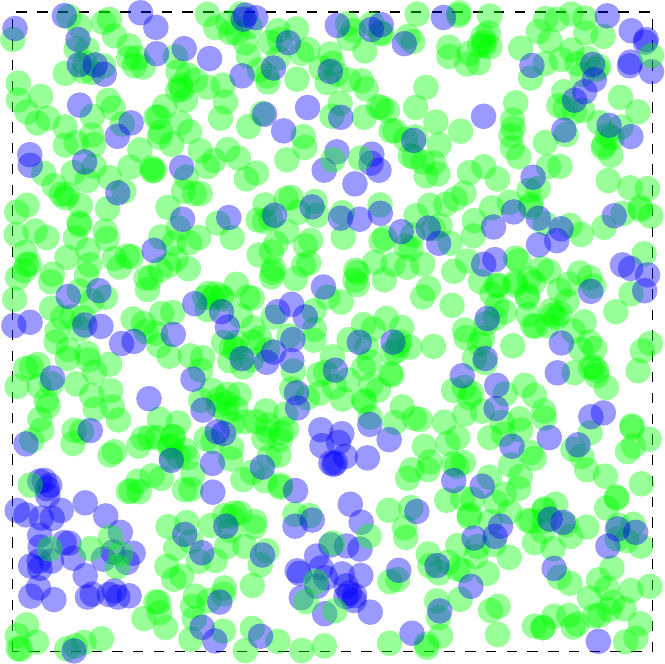}\hspace{1cm} }
 \end{center}
 We will apply the same stochastic dynamics also on the spatially discrete 
model on the lattice, and also to the mean-field model, at each finite system 
size $N$. 

We want to understand better the structural Gibbsian properties of the measures 
along a trajectory given by the time-evolution. For the purpose of concreteness we focus 
on the Euclidean model. 
We say the model shows a {\em dynamical Gibbs-non Gibbs 
transition} 
if the initial measure $\mu$ is Gibbs for a quasilocal specification, 
and for some time $t$ the time-evolved measure $\mu_t=\mu P_t$
 is not compatible with any quasilocal specification.
Here $P_t$ is the semigroup giving the distribution to find an infinite-volume configuration after time $t$ when starting with a given initial configuration, 
which is integrated over with respect to the starting measure $\mu$. 
In our example above $P_t$ is the symmetric independent spin-flip dynamics, and 
does not involve a randomization of the spatial degrees of freedom. However  
it is clear that one would like to study more generally also 
dependent dynamics, and also possibly irreversible dynamics, compare 
\cite{Li85,JaKu15}. 
Such studies have been performed at first for the Ising model on the lattice, for work 
on this and related work see \cite{EnFeHoRe02,KuRe06, KuLe07, EnKuOpRu08, 
ErKu10,EnFeHoRe10,FeHoMa14,
KrReZu17}.   

\subsection{Relation  to disordered systems}
To fix ideas, let us go to the lattice setup. 
{\em That} a time-evolved lattice measure $\mu_t$ 
(in our  the lattice Widom-Rowlinson model 
under symmetric spin-flip) is non-Gibbs is 
indicated by very long-range dependencies in its {\em conditional probabilities}, 
that is  \begin{equation}\begin{split}\label{acht}
\eta\mapsto \mu_t(\eta_i|\eta_{\Z^d\ba i})\end{split}
 \end{equation}
  behaves discontinuously w.r.t. the local topology.  
More precisely, as the r.h.s. is only defined up to measure-zero sets, this means 
that there is no  {\em version} which is continuous.

A useful strategy (at least for independent dynamics) is the following.  
Consider the {\em two-layer measure,} that is the joint distribution of spins and time zero 
and at in the future at time $t>0$ given by 
\begin{equation}\begin{split}\label{neun}\bar \mu_t(d\o,d\eta)=\mu(d\o)P_t(\o,d\eta).  \end{split}
 \end{equation}

Analyze {\em hidden phase transitions}
in {\em first-layer measure constrained on the future configuration} $\eta$. By this we mean the measure 
\begin{equation}\begin{split}\label{zehn}\bar \mu_t(d\o|\eta).\end{split}
 \end{equation}
 
A relation to disordered systems with quenched order appears when 
we view the configuration  $\eta=(\eta_i)_{i\in \Z^d}$ in the role of a  quenched disorder configuration.  By playing with (suitable finite-volume approximations of) 
conditional probabilities a picture emerges in which 
absence of phase transitions in the first-layer model implies Gibbsian behavior 
of the time-evolved model.  
The opposite implication, that the presence of a phase transition in the first-layer model 
implies the absence of the Gibbs property of the time-evolved model is true in many examples, 
and proved for a specific class of mean-field systems. 

Methods which are different from the two-layer picture 
are used for dependent dynamics. For mean-field systems and Kac-systems 
there are also path-large-deviation principles available which lead to fixed-end-point 
variational problems for trajectories of empirical measures. While in an abstract sense 
this is a solution, the analytical understanding of the structure of minimizers of such 
problems can be quite hard (see however \cite{HeKrKu17}). 
It is an open challenge 
to fully develop the analogous theory on the lattice, with ideas as suggested 
in \cite{EnFeHoRe10}.

\subsection{Results on dynamical GnG transitions for Euclidean model}

A marked infinite-volume configuration $\oo\in\OO$ (that is a signed point cloud) 
is called {\em good} for a 
specification $\g$ iff for any Euclidean ball $B$ we have 
\begin{equation*} 
\bigl| \g_{B}(f |\oo'_{B^{\rm c}})-\g_{B}(f |\oo_{B^{\rm c}})\bigr| \to0
\end{equation*}
as $\oo'\Rightarrow\oo$ in the sense of local convergence. 

We denote by 
$\OO(\g)$ the set of  good configurations for the specification $\g$. 
We say that a specification 
$\g$ is called {\em quasilocal} iff $\OO(\g)=\OO$, that is if all point clouds are good. 

A measure on signed point clouds 
$\mu$  is called {\em ql } (quasilocally Gibbs)  iff there exists a specification 
$\g$ such that $\OO(\g)=\OO$. 
$\mu$ is called {\em asql} (almost surely quasilocally Gibbs)  iff there exists $\g$ for which 
at least   $\mu(\OO(\g))=1$, that is the good points are a full-measure set.

Let us now 
describe the results on {\em Gibbsian transitions in time and intensity for $\mu^+$} obtained 
in \cite{JaKu17}. The measure $\mu^+$ is the measure in the Euclidean Widom-Rowlinson 
model obtained as an infinite-volume limit of finite-volume measures
with 
 the maximal boundary condition of overlapping plus-discs. 
By FKG (stochastic monotonicity) arguments this measure exists for all choices of parameters. 
We define
\begin{equation}\begin{split}\label{elf}t_G:=\frac{1}{2}\log\frac{\l_++\l_-}{\l_+-\l_-}\end{split}
 \end{equation}
for $\l_+>\l_-$. It will serve as a {\em reentrance time} into the Gibbsian region. 
We say that the model with intensity parameters $\l_+,\l_-$ is in the 
{\em high-intensity (percolating) regime} iff  
$\mu^+(B\leftrightarrow\infty)>0$ for some ball $B$  
(there is positive probability that there is an infinite cluster of overlapping discs containing $B$). 
Then the behavior of the time-evolved measure $\mu^+_t$ is summarized in the following table. 
\vskip0.5cm
{\centering
\begin{tabular}{c|ccc|c}
\hline
    & \multicolumn{3}{c|}{$\l_+>\l_-$}  & $\l_+=\l_-$ \\
 time   &$0<t< t_G$  & $t=t_G$ & $t_G<t\le\infty$
    & $0<t\le\infty$
 \\
 \hline
 high intensity & non-asql & asql, non-ql & ql &  non-asql
 \\
 low intensity & asq, non-ql & asql, non-ql & ql &  asql, non-ql
 \\
 \hline
 \end{tabular}
 \vskip0.5cm}
Main striking features are the 
immediate loss of the Gibbs property, and the appearance of  
full-measure sets of bad configurations (discontinuity points of any specification).  
More precise statements and detailed proofs can be found in \cite{JaKu17}.  
The proofs use  cluster representations for conditional probabilities of the time-evolved measure. 

The appearance of typical bad configurations can be heuristically understood: Infinite clusters 
in the time-zero model, together with the requirement that overlapping disks have 
the same sign, gives a strong rigidity in the first-layer model constrained on the future configuration 
$\eta$ at time $t$. Indeed, conditional 
on fixed locations in the percolating cluster, this cluster  
can only carry uniform plus signs, or unifom minus signs, at time zero.  
Keeping locations in a conditioning $\eta$ 
fixed and varying the signs arbitrarily far away provides then a very effective 
mechanism to induce a phase transition 
in the first-layer model. One shows that this implies jumps in certain conditional 
probabilities at time $t$. Hence every percolating point cloud can be a bad configuration, in  
the appropriate parameter regimes. As percolation is typical, this implies full-measure sets 
of bad configurations. 

\subsection{Results on dynamical GnG transitions for the mean-field Widom-Rowlinson model}
 Before we come to dynamics, we need to describe the equilibrium behavior of the mean-field model. 
The pressure of the mean-field model can be computed using large-deviation techniques 
(Varadhan's lemma) in terms of a variational formula where extremal points need to 
be found in the space of single-site probability distributions. Denoting by 
$L^1_{N}$ the fraction of spins with spin $1$ at system size $N$, and using similar notations 
for the other spin values,  we have 
\begin{equation*}
\begin{split}
		p(\beta,\alpha)&= \lim_{N\rightarrow \infty} \frac{1}{N} \log\int  e^{-N \beta L^1_{N}(\omega)L^{-1}_{N}(\omega) 
		} \prod_{j=1}^N \alpha(d\omega_j)\cr
		&= \sup_{\nu\in \mathcal{M}_1(\{-1,0,1\})} (- \b \nu(1)\nu(-1)- I(\nu\vert \alpha ))
\end{split}
\end{equation*}	
where $I(\nu\vert \alpha )$ is the relative entropy of a single-site distribution $\nu$ w.r.t. 
the a priori measure $\a$. Correspondingly, the empirical distribution satisfies a large deviation 
principle with speed $N$, and rate function given 
by the negative of the expression under the $\sup$ plus a suitable constant. 

A discussion of the variational problem (see \cite{KiKu18}) shows:  
The symmetric model at any $\alpha(1)=\alpha(-1)>0$
has a second-order phase transition driven by repulsion strength $\beta>0$ 
at critical repulsion strength given by $\beta_c= 2+e \frac{\alpha(0)}{\alpha(1)}$. 
%


An explicit solution for the mean-field Widom-Rowlinson model at time $t=0$ can be obtained as follows. 
We parametrize the empirical spin distribution $\nu$ via coordinates $(x,m)$ where  
$x$ plays the role of occupation density, and $m$ plays the role of the 
magnetization on the occupied sites, writing 
\begin{equation} \bigl(\nu(-1), \nu(0), \nu(1) \bigr)=( 
 \frac{x}{2}(1-m),1-x, \frac{x}{2}(1+m))
\end{equation}	
%
Next we parametrize the a priori measure $\alpha$ via coordinates $(h,l)$
where $h= \frac{1}{2} \log\left(\frac{\alpha(1)}{\alpha(-1)}\right)$  plays the role 
of an external magnetic field, and 
 $l:=\log\frac{1-\alpha(0)}{\alpha(0)}$   describes a bias on the occupation variables. 
Using these coordinates, the pressure can be written as 
\begin{equation*}\begin{split}
 &p(\beta,\alpha)=\log\alpha(0)+ \sup_{0\leq x,\vert m\vert\leq 1} \Bigl( \underbrace{-\frac{\beta x^2}{4}+x(l-\log(2\cosh(h)) -J(x)}_{\hbox{\small
 part for occupation density}}
 \\
& + x (\underbrace{\frac{\beta x m^2 }{4} +h m-I(m)-\log 2 }_{\hbox{\small Ising part at occupation-dependent temperature}}) \Bigr)
 \end{split}
 \end{equation*}
  with an entropy for spins 
 $ I(m)= \frac{1-m}{2}\log(\frac{1-m}{2})+\frac{1+m}{2}\log(\frac{1+m}{2}) $ 
  and an entropy for occupation variables 
  $J(x) = (1-x)\log(3(1-x))+ x\log(\frac{3x}{2})$. 
It turns out that the symmetric antiferromagnetic model ($\b<0$) has a first-order transition when crossing a line in $\b,\alpha(0)$-space, where  
jumps occur in occupation density $x$, at fixed zero magnetization $m=0$. 

%

We are mostly interested in the ferromagnetic model. 
In this situation we obtain that repulsion parameter $\b>0$,  
a priori measure $\alpha=\alpha(h,l)$, 
and typical values $(m,x)$ of the empirical distribution, 
 are related via the parametrization 
 \begin{align*}\label{equi: asym beta}
\beta=\beta(m;\alpha)&= \frac{2}{m}(I'(m)-h)  (1+e^{-l+\log(2\cosh(h)) +\frac{1}{m}(I'(m)-h)-mI'(m)+I(m)})\cr
x=x(m;\alpha)&=(1+e^{-l+\log(2\cosh(h)) +\frac{1}{m}(I'(m)-h)-mI'(m)+I(m)})^{-1} 
 \end{align*}
We remark as a corollary that the model has {\em mean-field critical exponents:} 
Fix any $\alpha(0)\in (0,1)$: 
Let $\beta_c$ be the corresponding critical value for the symmetric model. Then there are positive constants 
such that   
\begin{equation*}\label{equi: crit exp beta}
 	\lim_{\beta\downarrow \beta_c}\frac{m(\b,h=0)}{(\beta-\beta_c)^\frac{1}{2}}=
 	c , \qquad	\lim_{h\downarrow 0}\frac{m(\b_c,h)}{h^\frac{1}{3}}= c'
 		\end{equation*}
The main point is the study of the {\em dynamical Gibbs-non Gibbs transitions} 
under rate-one symmetric independent spin-flip, keeping holes fixed, 
according to transition probabilities (\ref{sieben}). 
 
Recall the notion of sequentially Gibbs and the notion of bad empirical measure, see 
\eqref{fuenfa},\eqref{sechs} and the text below. For the sake of this review 
let us just present the time evolution of the set of bad empirical measures in the regime of an inverse 
temperature of the time-zero model in the region of 
$\b>3$ in a plot (compare \cite{KiKu18} for the full statement 
of the theorem describing all dynamical transitions).  In the plot the inverse temperature 
of the time zero model is $\b=5$ and we are 
starting from a symmetric model with $\a(+)=\a(-)$. 

Here time increases from the top left to the top right, then from bottom left to bottom right.  
The main features are the following:  
There is a  short-time Gibbs regime for all $\beta,\alpha$.  
Small repulsion strength $\beta\leq 2$ implies the Gibbs property of the time-evolved model 
for all times $t$. 
The set of bad empirical measures at given $\b$ 
has dimension one as a subset of the simplex, in the interior of 
its existence time-interval. It 
can be a union of disconnected curves, a branching curve (which has a Y-shape, see picture), 
or a line which is growing with time  (growing antenna).  

\vskip0.5cm
\includegraphics[trim = 0mm 40mm 0mm 0mm  , clip,width=0.3\textwidth]{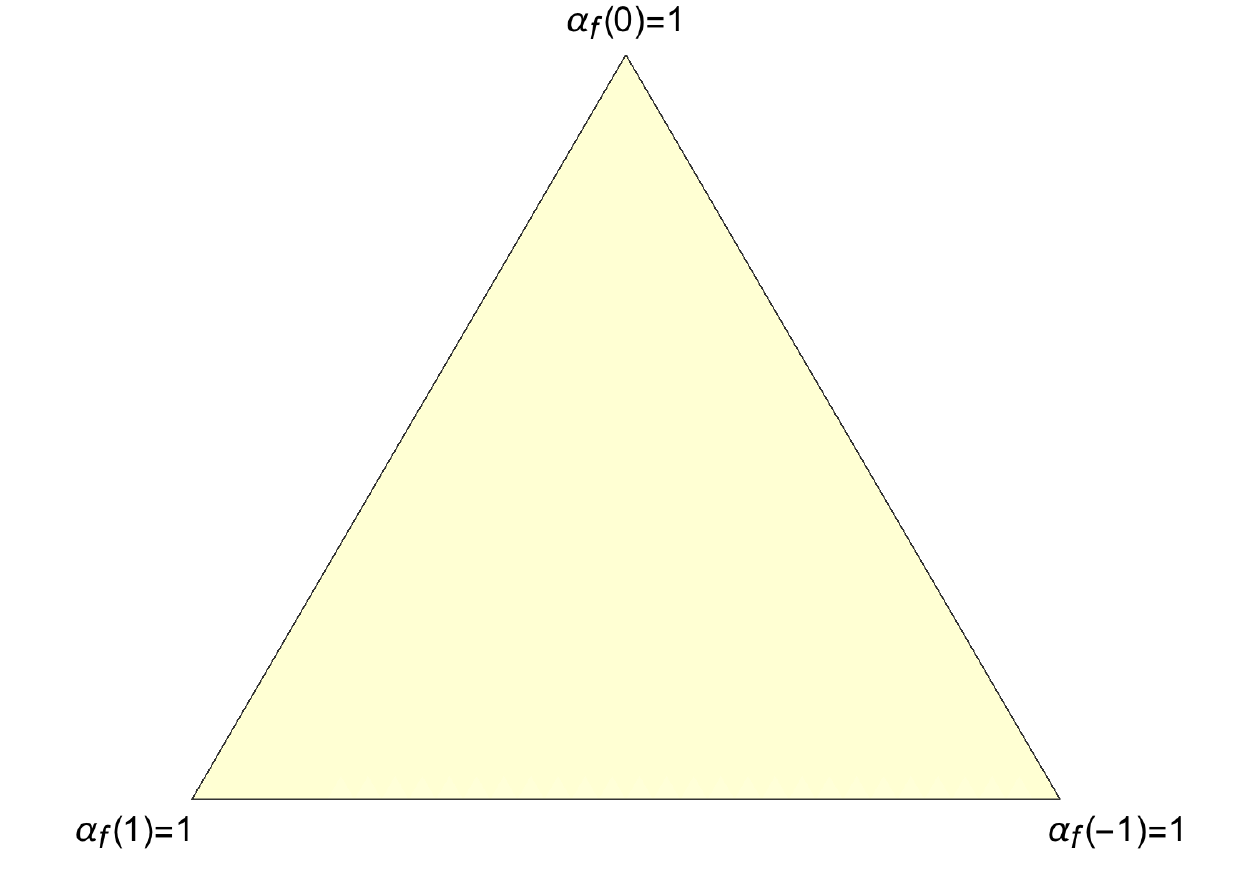}
 \hfill
\includegraphics[trim = 0mm 40mm 0mm 0mm  , clip,width=0.3\textwidth]{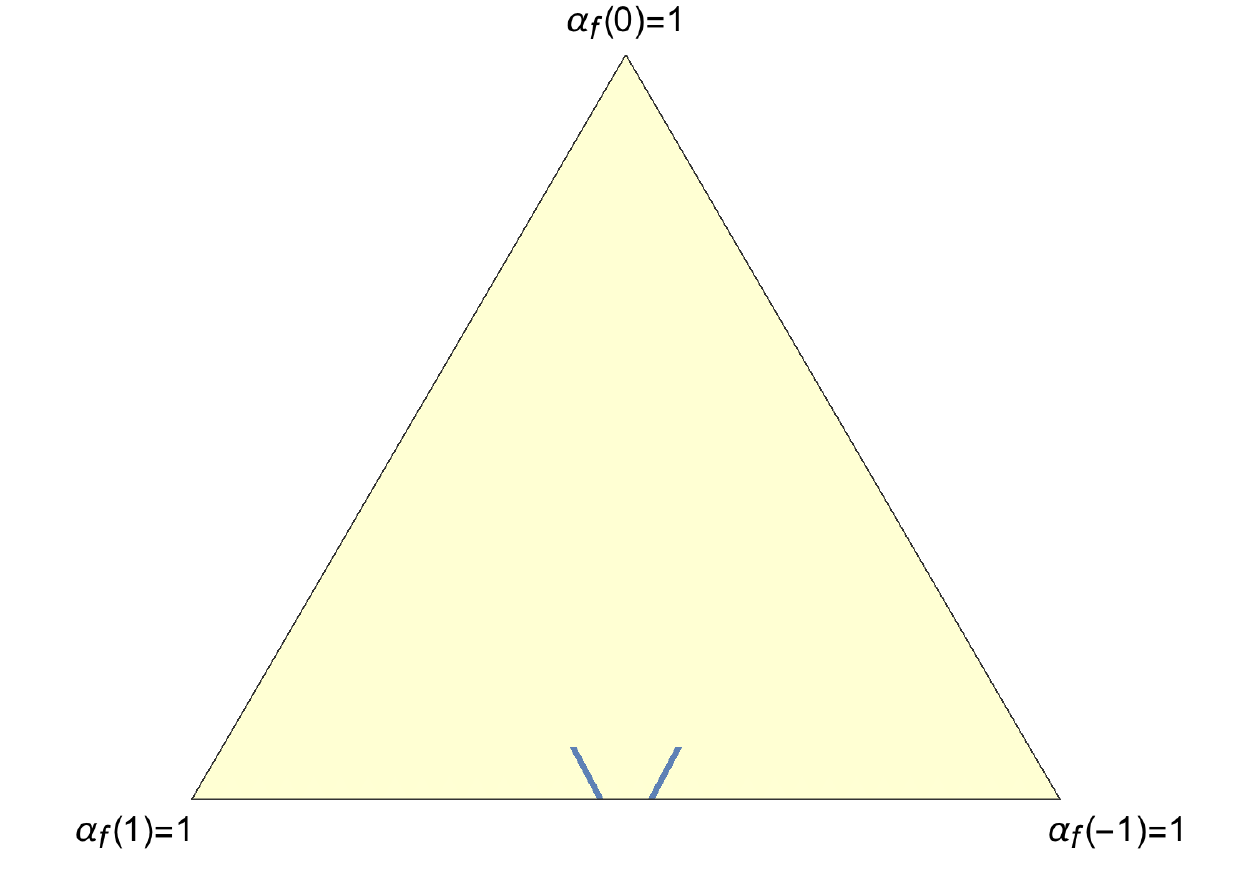}
\hfill 
\includegraphics[trim = 0mm 40mm 0mm 0mm  , clip,width=0.3\textwidth]{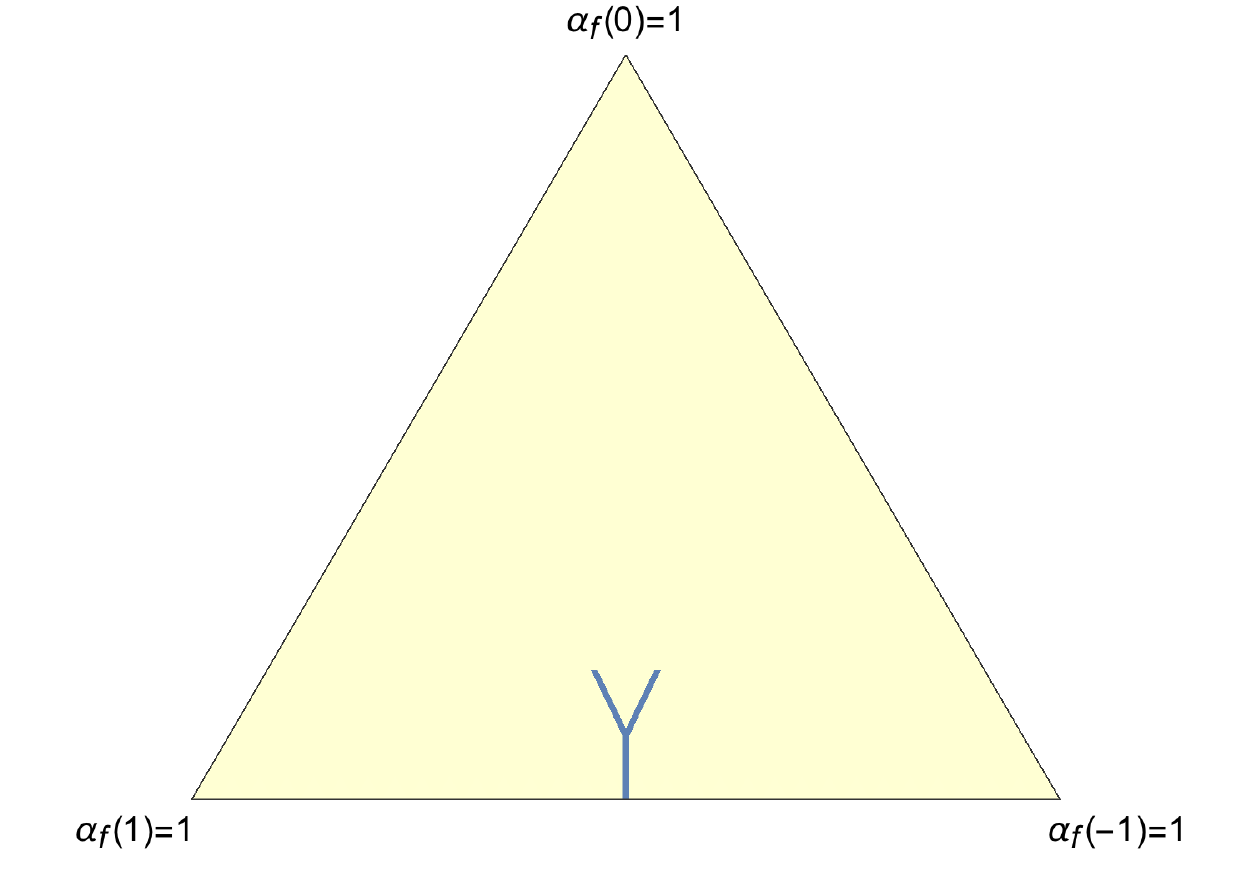}

\vskip0.5cm
 \includegraphics[trim = 0mm 40mm 0mm 0mm  , clip,width=0.3\textwidth]{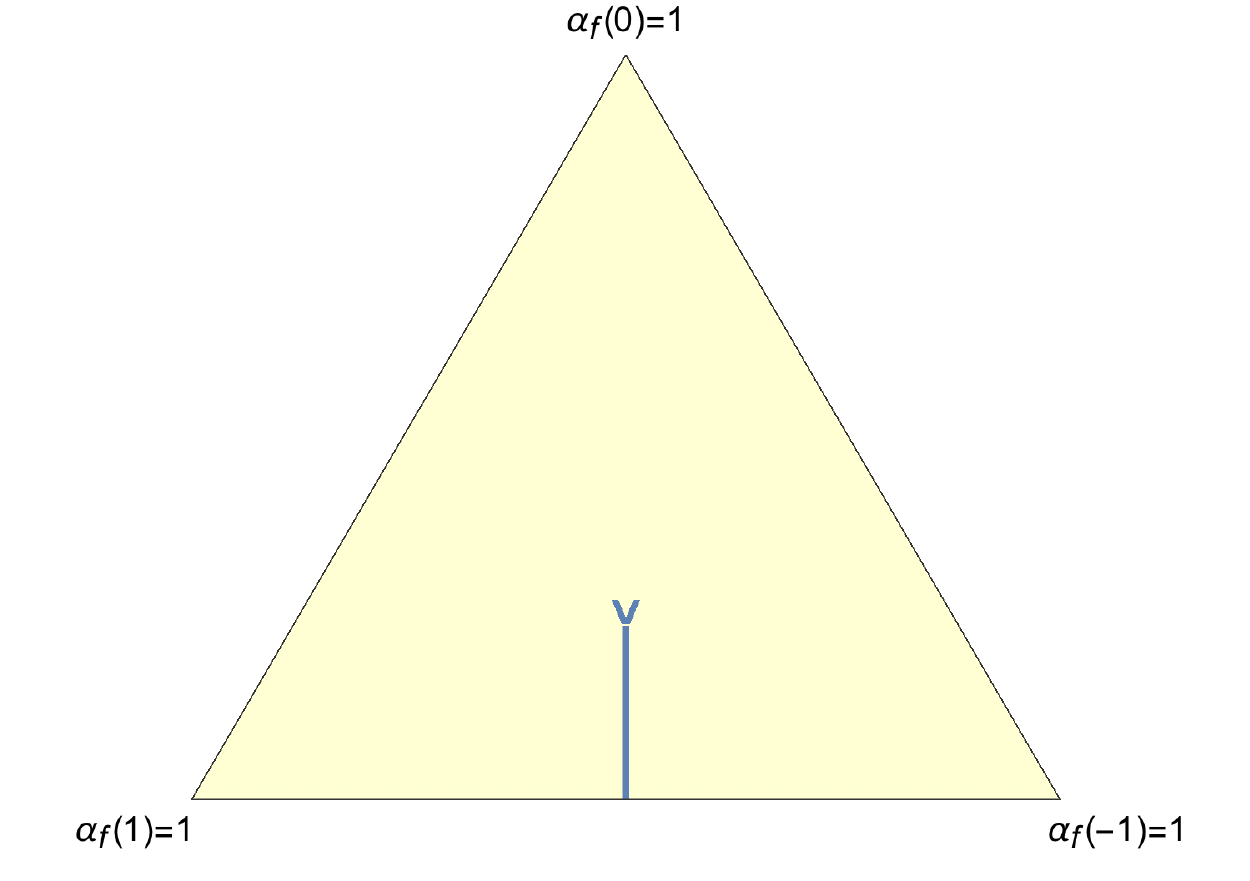}
 \hfill 
\includegraphics[trim = 0mm 40mm 0mm 0mm  , clip,width=0.3\textwidth]{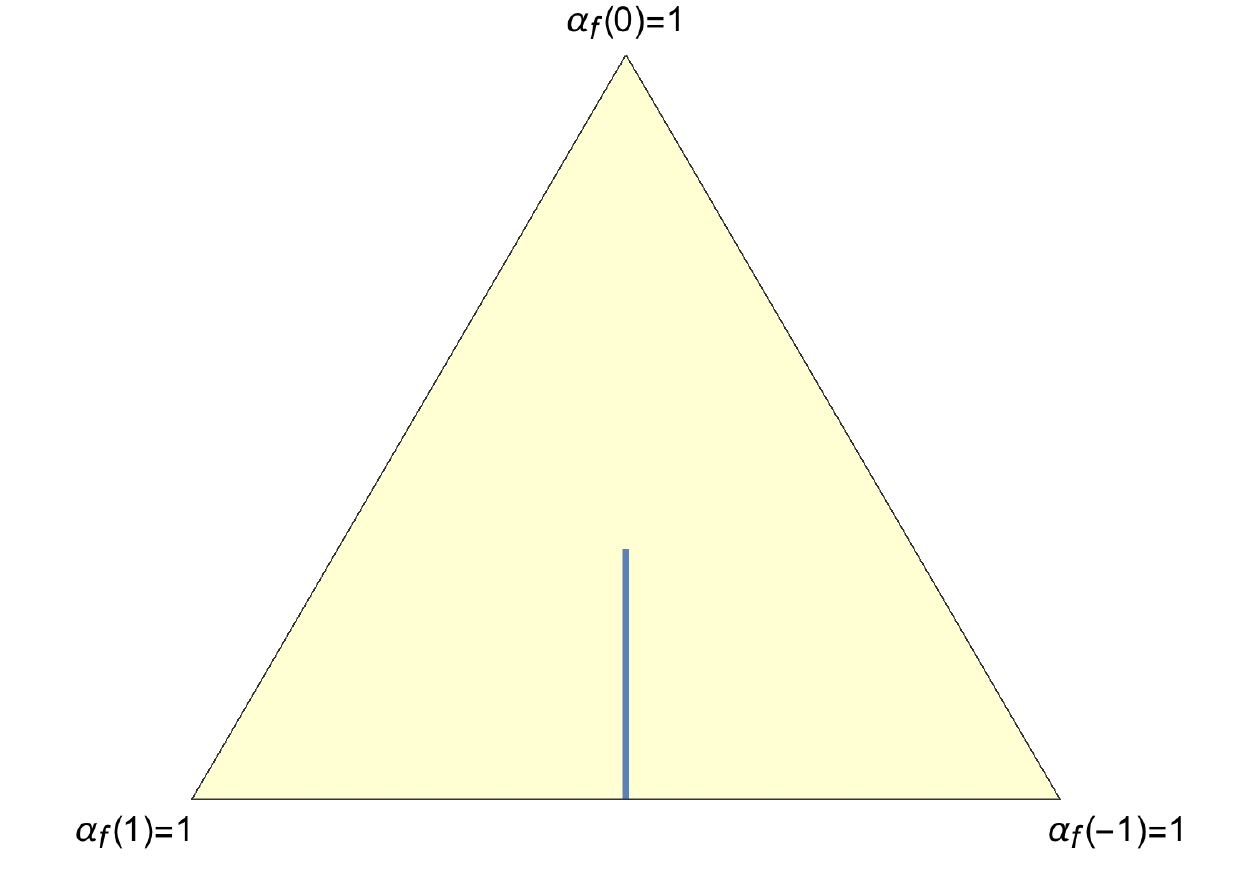}
\hfill 
\includegraphics[trim = 0mm 40mm 0mm 0mm  , clip,width=0.3\textwidth]{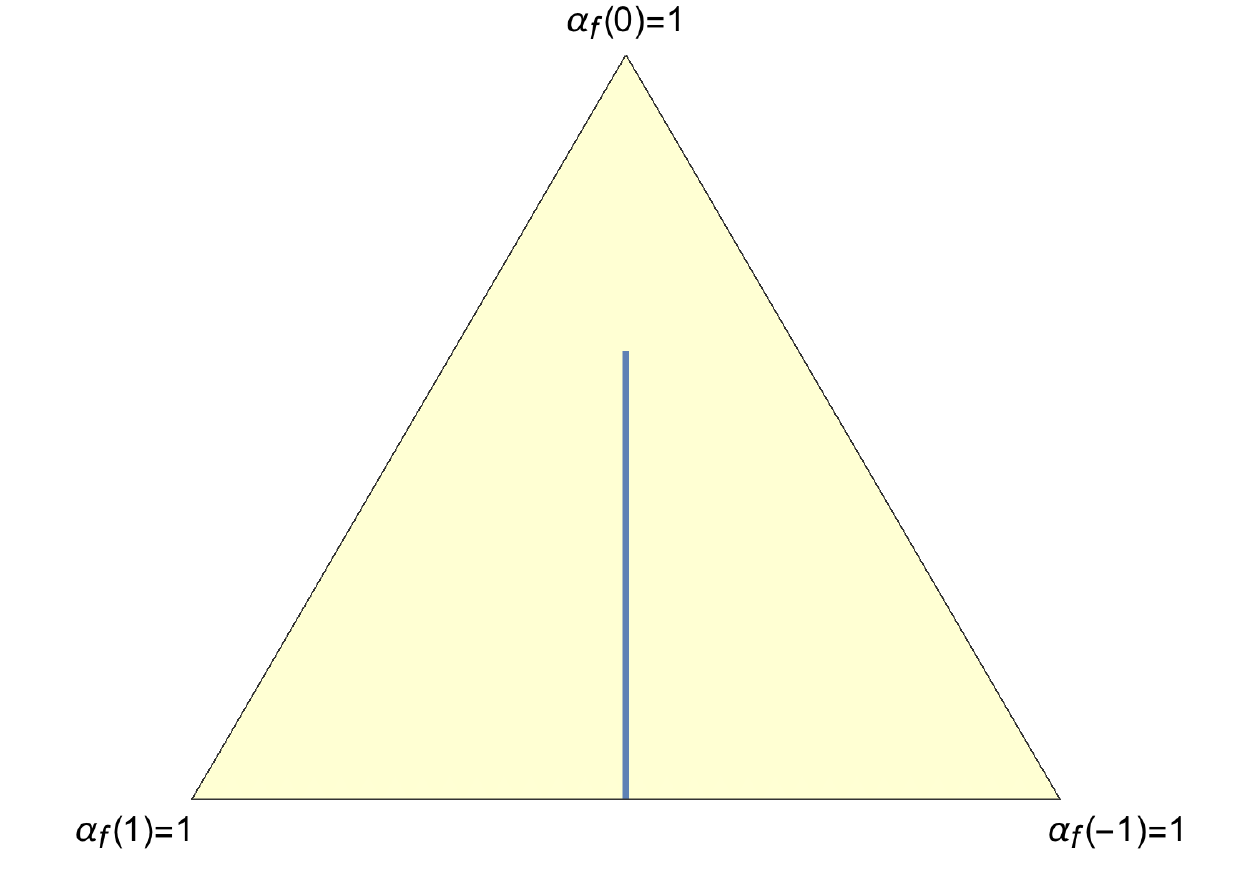}

These features can be understood by a bifurcation analysis of a rate function 
describing the first-layer model constrained on an empirical distribution. 
For our model, this analysis is closely related to that of a time-evolved Ising model in 
the following way. 
The bad measures $\alpha_f$ in the time-evolved
mean-field WiRo model after time $t$ satisfy for any symmetric a priori measure $\alpha$
 \begin{equation*}
  \begin{split}
  &B^{WiRo}(\beta,t)=\Bigl\{\alpha_f\in \mathcal{M}_1(\{-1,0,1\}), \,\,
   \frac{\alpha_f(1)-\alpha_f(-1)}{\alpha_f(\{1,-1\})}\in 
  B^{Ising}\bigl(\frac{\beta \alpha_f(\{1,-1\}\bigr)}{2},t)
  \Bigr\}.
  \end{split}
  \end{equation*}
where $ B^{Ising}(\beta_I,t)$ denotes the set of 
bad magnetizations for the time-evolved Curie-Weiss Ising model 
with initial inverse temperature $\b_I$.   It is known from \cite{KuLe07} that 
$ B^{Ising}(\beta_I,t)$ turns out to be either empty, contain the magnetization value zero, 
or to be given by a symmetric pair. 


What can we say about the {\em typicality} of bad points in the time-evolved mean-field 
Widom-Rowlinson model? Typicality means in the mean-field context 
that the minimizers of the large-deviation 
rate function of the time-evolved model are contained in the set of bad magnetizations. 
It is an  analytical principle for time-evolved mean-field Ising models  that there is an 
atypicality of bad configurations. This follows from the principle of 
{\em preservation of semiconcavity} for time-evolved rate functions which are defined
 via integrals over Lagrange densities (\cite{KrReZu17}). In simple words this regularity 
 statement means that kinks in a rate function can never appear at local minima. 

Our model does not fall in the Ising class, but the corresponding statement can be proved 
explicitly. It is very nicely illustrated in our model by the following plot. The repulsion strength 
of the model at time zero is $\beta=4>3$ is the low-temperature region. 
The red Y-shaped set denotes the {\em set of bad empirical measures} at a fixed intermediate 
time. Its form is independent of the initial a priori distribution $\a$, as long as this is symmetric.  
By comparison the {\em typical configurations} for any $\a$, after time $t$ are solid blue. 
They arise as time-evolution of the dotted blue lines describing typical measures at time zero. 

\begin{center}
\includegraphics[trim = 0mm 40mm 0mm 0mm  , clip,width=0.45\textwidth]{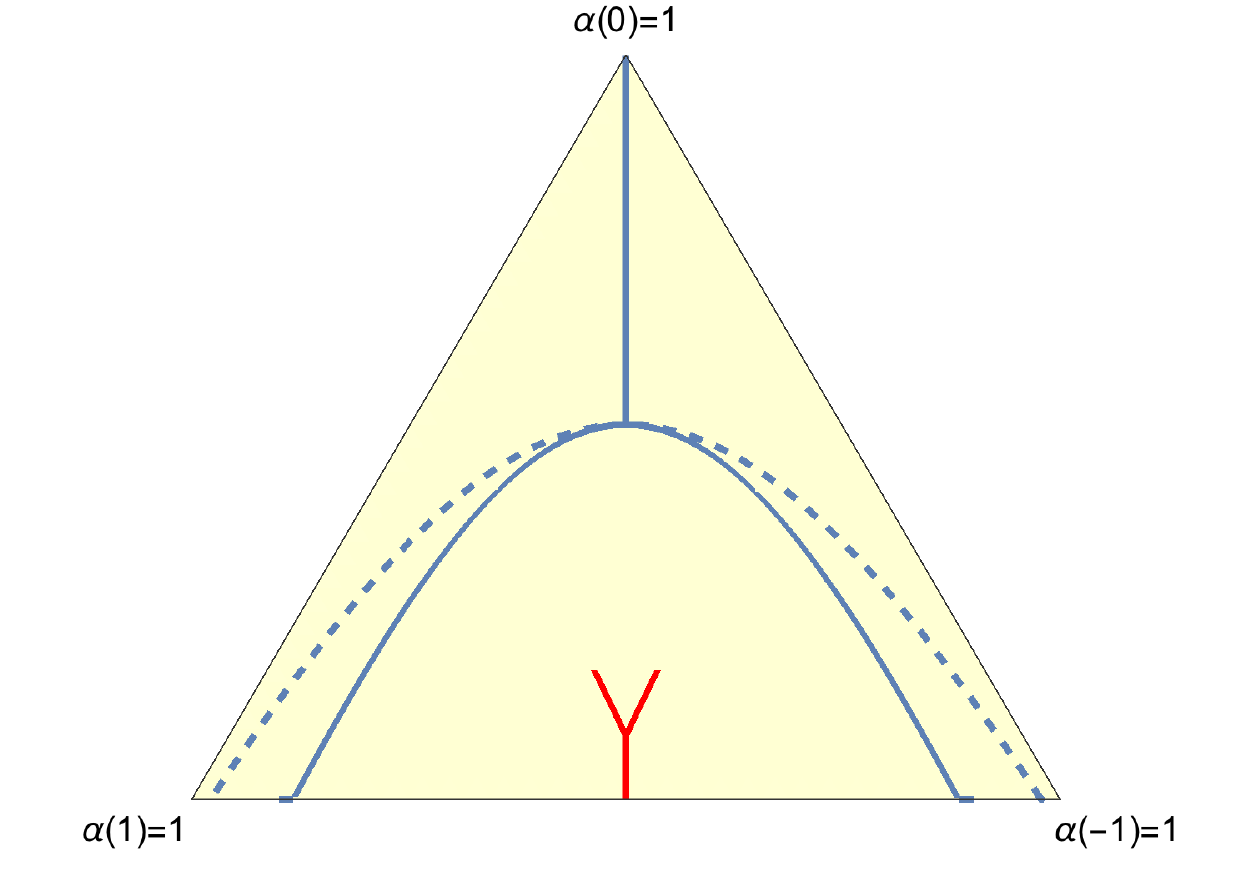}
\end{center}

\subsection{Lattice Widom-Rowlinson model under time-evolution}

We describe mainly the {\em lattice soft-core Widom-Rowlinson model}. To prove 
the Gibbs property of the time-evolved model in appropriate regions, 
Dobrushin-methods are useful, as we outline 
now. 
Let $\gamma:=(\gamma_{\Lambda})_{\Lambda\Subset \Z^d}$ be a quasilocal specification
on the lattice. The Dobrushin interdependence matrix, is defined by  
	\begin{align}
		C_{ij}(\gamma) = \sup_{\omega_{\Z^d\backslash \{j\}}=\eta_{\Z^d\backslash \{j\}} }\Vert \gamma_{\{i\}}(\cdot\vert\omega)-\gamma_{\{i\}}(\cdot\vert\eta)\Vert_{TV,i}
	\end{align}
for sites $i\neq j$. 	
The main theorem, due to Dobrushin, states: 
If the {\em Dobrushin condition} holds, namely if 
$c(\gamma):= \sup_{i\in\Z^d} \sum_{j\in \Z^d} C_{ij}(\gamma)<1$,  
then $\vert \mathcal{G}(\gamma)  \vert = 1$. The theory also allows to control 
the unique Gibbs measure under perturbations of the specification, understand correlation
decay in the measure, and derive more useful consequences \cite{Ge11}. 
We show the {\em Dobrushin region } (the region in parameter space for which $c(\gamma)<1$) 
for the spatially homogeneous soft-core model on $\Z^2$.  The plot shows, for  different values 
of the repulsion strength $\beta$, the 
Dobrushin region (dark shaded) in the space of a priori measures $\alpha\in \mathcal{M}(\{-1,0,1\})$, 
projected to the $\a(1),\a(-1)$-plane. 
 \vspace{0.5cm}
 
\hspace{1cm} 
	\includegraphics[width=70pt,height=70pt]{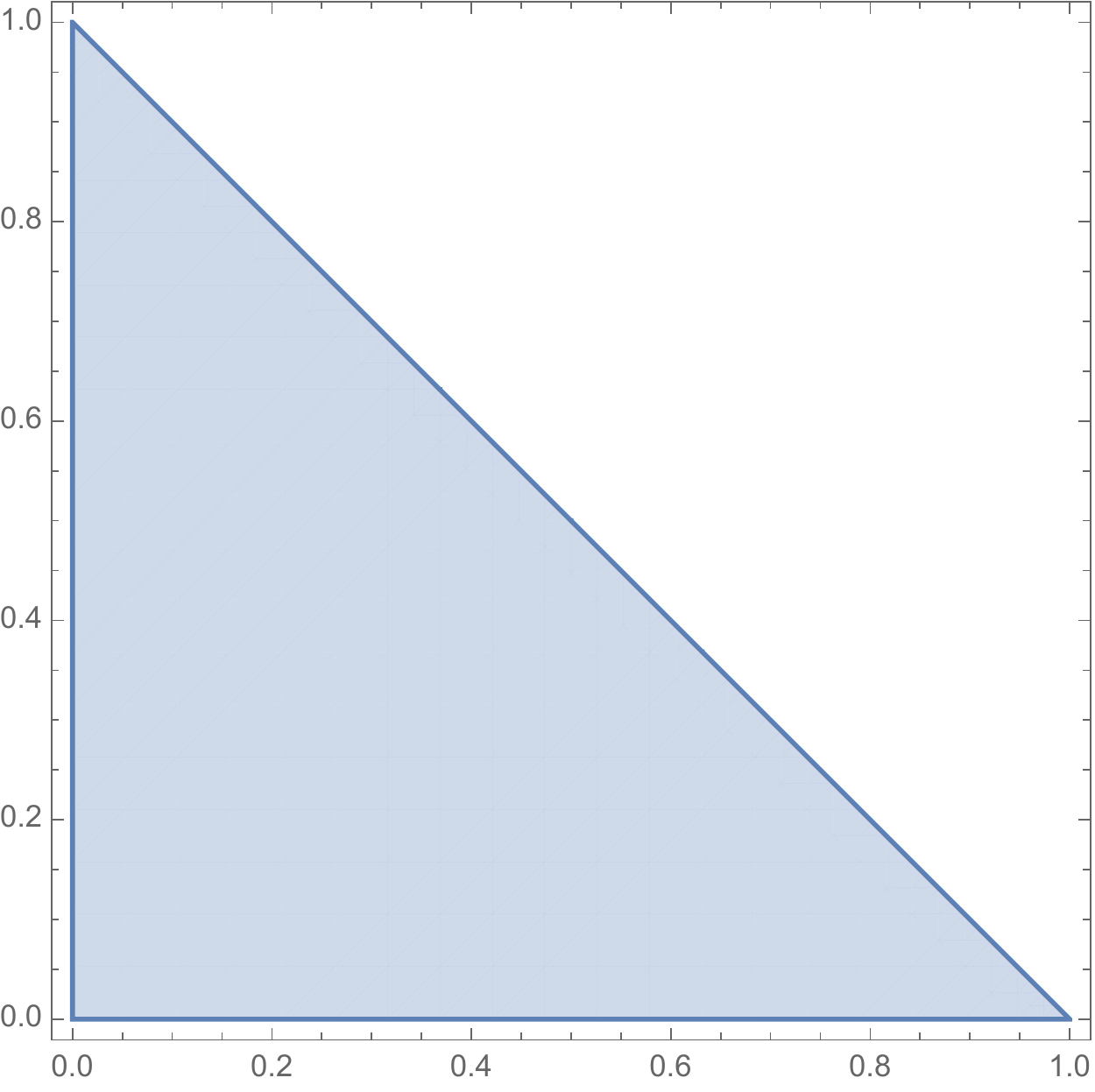}
\hspace{1cm}  
\includegraphics[width=70pt,height=70pt]{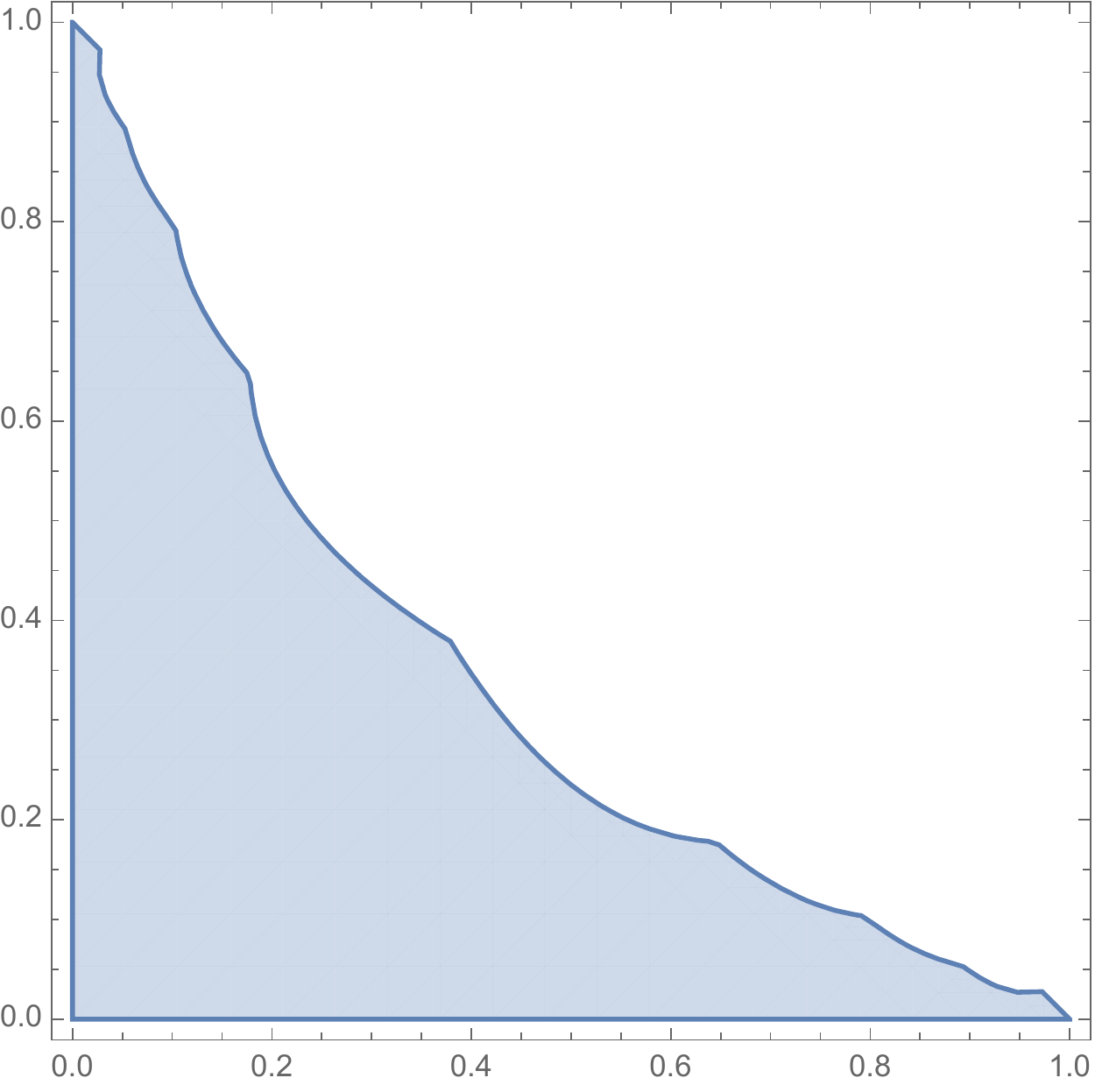}
\hspace{1cm} \includegraphics[width=70pt,height=70pt]{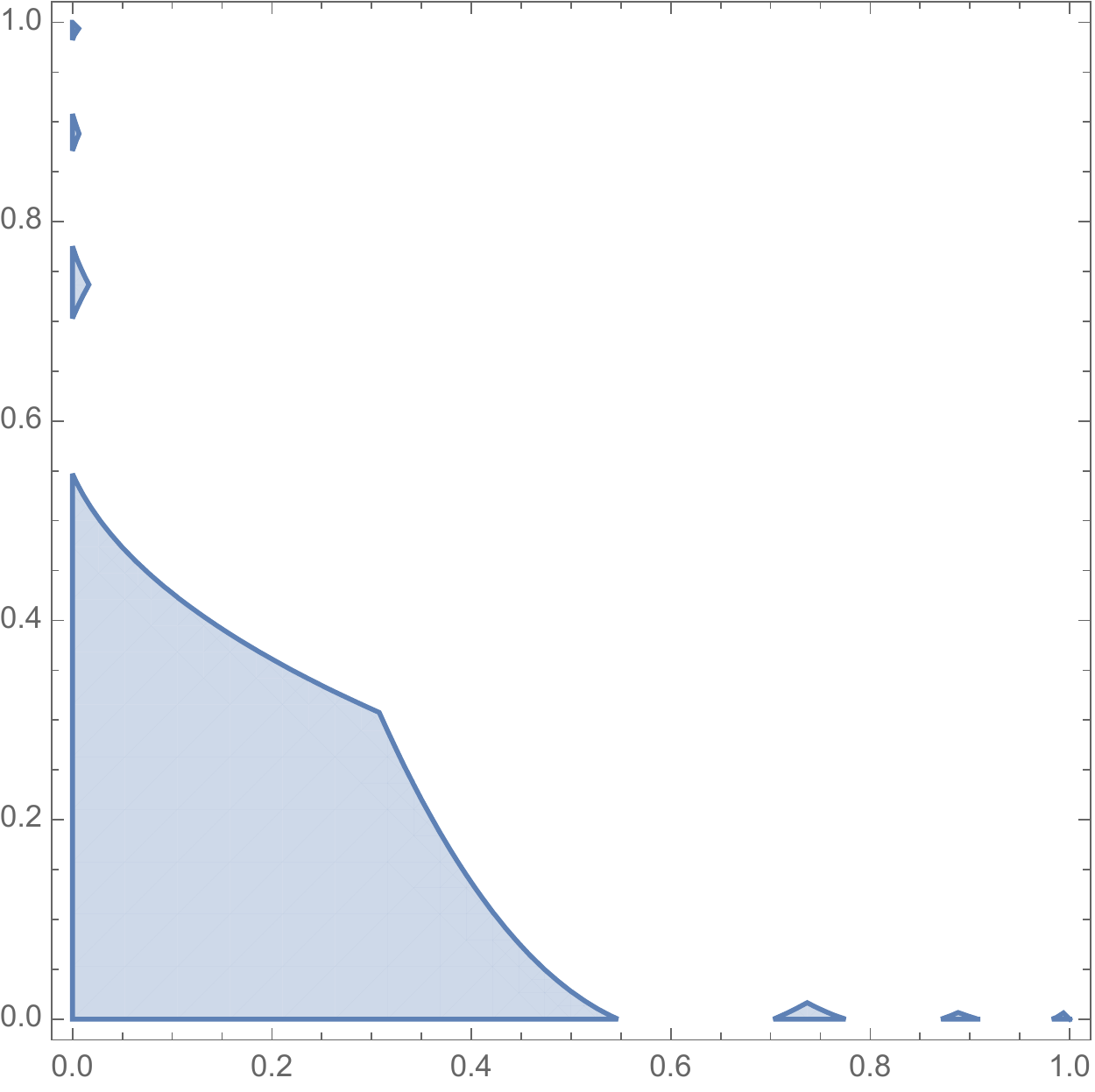}
	\hspace{1cm} 
	\includegraphics[width=70pt,height=70pt]{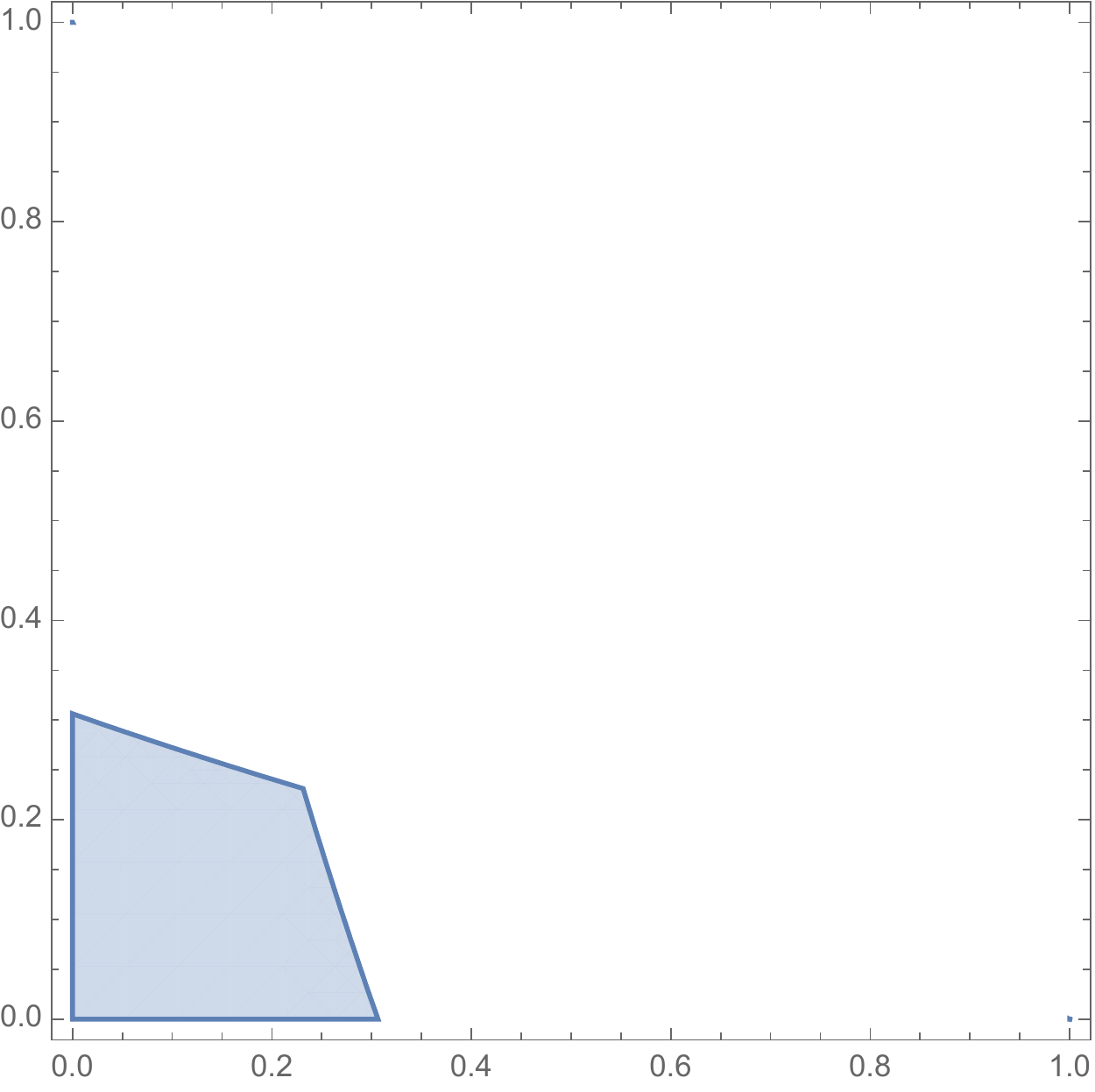}\\
$\phantom{1111111111}\beta<\frac{1}{2}\phantom{11111111111111}\beta=0.75\phantom{111111111111}\beta=1.05\phantom{111111111111}\beta=2$ 
\vspace{0.5cm}

It turns out that the boundary of the Dobrushin region consists of finitely many pieces 
of ellipses. We have the following theorem for the homogeneous model: \\
(1) Let $0\leq \beta d<1$. Then for all $\alpha\in\mathcal{M}(\{-1,0,1\})$ the soft-core Widom-Rowlinson model satisfies the Dobrushin condition.\\
 (2)	For every $\beta>0$ there exists an $\epsilon:=\epsilon(\beta)>0$ such that the soft-core model satisfies the Dobrushin condition if  $d_{TV}(\alpha,\delta_1)<\epsilon$ or $d_{TV}(\alpha,\delta_{-1})<\epsilon$.\\
 With Dobrushin techniques one  controls not only the translation-invariant model, 
 but also the first-layer model 
constrained on future configurations. With this one may prove also the {\em Gibbs property for the time-evolved 
model for small times}: 
		Let $\alpha \in \mathcal{M}(\{-1,0,1\})$, $ \beta\geq 0$,  
		and  let $\mu\in\mathcal{G}(\gamma_{\beta,\alpha})$ 
		by any Gibbs measure. Then there exists a time $t_c>0$ such that for all $t<t_c$ the time-evolved measure $\mu_t$ is a Gibbs measure for some quasilocal specification $\gamma_t$. For the proof see \cite{KiKu19}, it extends methods of 
		\cite{KuOp08} to a situation of degenerate time evolutions (where not all transitions 
		are allowed)  to control {\em all} first-layer models for possible to control all first-layer models 
		\eqref{zehn}  for possible end-conditionings $\eta$. 
		We remark that the method 
does not make use of the lattice structure, but gives the same result of short-time Gibbsianness for any graph with bounded degree, for instance a regular tree.  
	
For the opposite direction we prove: In the soft-core model on the lattice, at sufficiently large 
repulsion strength, the maximal measure $\mu^+_t$ is non-Gibbs,  for times $t$ which are sufficiently large. For the proof it suffices to exhibit one non-removable bad configuration for the 
single-site probability of the time-evolved measure. We may choose in our case 
a fully occupied checkerboard configuration of alternating plus- and minus-spins, 
and show that this configuration is bad, noting that we are  reduced basically to an Ising situation for 
this conditioning.

For the hard-core lattice model under time-evolution, Dobrushin techniques can not 
be applied, as some entries of the Dobrushin matrix will necessarily become equal to one. 
This is not just a shortcoming of the proof. Indeed, we find an immediate loss of 
the Gibbs property, as in the Euclidean model, for the proof see \cite{KiKu19}. 

\subsection{The Widom-Rowlinson model on a Cayley tree}

Let us now for our graph consider a Cayley tree, which is the infinite graph which has no 
loops, and where each vertex has precisely $k+1$ nearest neighbors. The 
Widom-Rowlinson model in the hard-core version, and in the soft-core version, is again defined 
by the specification kernels of \eqref{vier} and \eqref{viera}. \\
We need to start with a good understanding of the Widom-Rowlinson model 
in equilibrium. The tree-automorphism-invariant Gibbs measures which are also 
tree-indexed Markov chains (also known as tree-invariant splitting Gibbs-measures) are uniquely 
described via boundary laws, which appear as solutions of a parameter-dependent 
two-dimensional fixed point equation (appearing as a tree recursion). 
As a general abstract fact, extremal Gibbs measures
for tree models  
are always splitting Gibbs measures, the opposite is in general not true \cite{Ge11}. For certain 
classes of  {\em hard-core} models on trees 
the characterizations of solutions can be found in \cite{Ro13}, at least for low enough degree 
of the tree.  
For the equilibrium states of the ferromagnetic {\em soft-core} model on the Cayley tree 
we find the following \cite{KiKuRo19}.   
In the {\em antiferromagnetic model}  with symmetric intensities 
there is a transition in the hole-density, somewhat 
similar to that in the mean-field model briefly described above. It  
can be very explicitly analyzed for any order $k$, with explicit transition lines 
in the interaction-intensity diagram. For the {\em ferromagnetic model with symmetric intensities}, 
for the trees with $3$ and $4$ nearest neighbors,  
the critical lines for the ferromagnetic phase transition are again explicit, with complete description of 
all tree-invariant splitting Gibbs measures.  For higher $k$, there are only bounds on critical curves, 
which we conjecture to be sharp, see \cite{KiKuRo19}.    

What about spin-flip time evolution of these measures on the tree? The Gibbsian behavior of a time-evolved 
model can be very different from the behavior in other geometries. 
For the Ising model on a Cayley-tree under independent stochastic spin-flip in \cite{EnErIaKu12} 
the following was proved: 
The set of bad measures may depend on the choice of the initial Gibbs measure of the 
time-evolved state. There can be multiple transition times in the model with 
zero external magnetic field, and full-measure sets of bad measures. 
For the time-evolved Widom-Rowlinson model on the Cayley tree, this is an open problem.

%


%

\section*{Acknowledgements} 
I am very grateful I met Anton, for all the discussions I had with him, and for all 
inspiration he gave, during my Ph.D., in later years, until today.  
I wish him many many more years, I am looking forward  to 
many more of his contributions, to mathematics and beyond!

\end{document}